\documentclass[a4paper,12pt]{article}
\usepackage{amsmath}
\usepackage{amssymb}
\usepackage{amsthm}

\theoremstyle{plain}
\newtheorem{theorem}{Theorem}[section]

\theoremstyle{definition}

\def\ppmod{\mkern-16mu \pmod}
\newcommand*{\Z}{\mathbb{Z}}

\newcommand*{\C}{\mathbb{C}}

\begin{document}

\title{Discrete Ramanujan-Fourier Transform of\\ Even Functions (mod $r$)}
\author{
Pentti Haukkanen\\
Department of Mathematics, Statistics and Philosophy\\
FI-33014 University of Tampere\\
Finland\\
mapehau@uta.fi\\
\\
{\sl Published in Indian J. Math. Math. Sci. 3  (2007),  no. 1, 75--80.}}
\date{2007}
\maketitle

{\bf Abstract.} An arithmetical function
$f$ is said to be even $\pmod r$ if
$
f(n)=f((n,r))
$
for all $n\in\Z^+$, where $(n, r)$ is the greatest common divisor of 
$n$ and $r$. We adopt a linear algebraic approach to  
show that the Discrete Fourier Transform of an even function (mod $r$) 
can be written in terms of Ramanujan's sum and may thus be referred to as 
the Discrete Ramanujan-Fourier Transform. 

\medskip


\medskip


\section{Introduction}

By an arithmetical function we mean 
a complex-valued function defined on 
the set of positive integers. 
For a positive integer $r$, an arithmetical function $f$ is said to be 
{\sl periodic} $\pmod r$ if 
$f(n+r)=f(n)$ for all $n\in\Z^+$. 
Every periodic function $f$ $\pmod r$ can be written uniquely as 
\begin{equation}\label{y:IDFT}
f(n)=r^{-1}\sum_{k=1}^{r}F_f(k)\epsilon_k(n),   
\end{equation}
where 
\begin{equation}\label{y:DFT}
F_f(k)=\sum_{n=1}^{r}f(n)\epsilon_k(-n)  
\end{equation}
and $\epsilon_k$ denotes the periodic function $\pmod r$ defined as
$$
\epsilon_k(n)=\exp (2\pi ikn/r). 
$$
The function $F_f$ in (\ref{y:DFT}) is referred to as 
the Discrete Fourier Transform (DFT)
of $f$, and (\ref{y:IDFT}) is 
the Inverse Discrete Fourier Transform (IDFT).

An arithmetical function
$f$ is said to be {\em even} $\pmod r$ if
$$
f(n)=f((n,r))
$$
for all $n\in\Z^+$, where $(n,r)$ is the greatest common divisor of 
$n$ and $r$. It is easy to see that every even function $\pmod r$ is 
periodic $\pmod r$. 
Ramanujan's sum $C(n, r)$ is defined as 
$$
C(n, r)=\sum_{k\ppmod r\atop (k, r)=1}\exp(2\pi ikn/r)
$$ 
and is an example of an even function $\pmod r$. 

In this paper we show that the DFT (\ref{y:DFT}) and IDFT (\ref{y:IDFT}) of 
an even function $f$ $\pmod r$ can be written in a concise form using 
Ramanujan's sum $C(n, r)$, see Section 3. We also review a proof of  
(\ref{y:IDFT}) and (\ref{y:DFT}) for periodic functions $\pmod r$, 
see Section 2, and review (\ref{y:IDFT}) and (\ref{y:DFT}) for the 
Cauchy product of periodic functions $\pmod r$, see Section 4. 
The Cauchy product of periodic functions $f$ and $g$ $\pmod r$ is defined as 
$$
(f\circ g)(n)
=\sum_{a+b\equiv n\ppmod r} f(a)g(b).
$$

The results of this paper may be considered to be known. 
They have not been presented
in exactly this form and we hope that this paper will provide a clear approach to 
the elementary theory of even functions $\pmod r$.  

The concept of an even function $\pmod r$ originates from Cohen \cite{C}  
and was further studied 
by Cohen in subsequent papers \cite{C1,C2,C3}. 
General accounts of even functions $\pmod r$ can be found 
in the books by McCarthy \cite{M} and Sivaramakrishnan \cite{Si}. 
For recent papers on even functions $\pmod r$ we refer to \cite{Sa,T}. 
Material on periodic functions $\pmod r$ can be found in  
the book by Apostol \cite{A}.

\section{Proof of (\ref{y:IDFT}) and (\ref{y:DFT})}

Let $P_r$ denote the set of all periodic arithmetical functions 
$\pmod r$. 
It is clear that $P_r$ is a complex vector space under the usual addition and 
scalar multiplication. 
In fact, $P_r$ is isomorphic to $\C^r$. 
Further, $P_r$ is a complex inner product space under 
the Euclidean inner product given as 
\begin{equation}\label{y:innerp}
\langle f, g\rangle = \sum_{n=1}^r f(n) \overline{g(n)}
=(f\overline{g}\circ\zeta)(r),
\end{equation}
where $\zeta$ is the constant function $1$. 
The set $\{r^{-1/2}\epsilon_k:\ k=1, 2,\ldots, r\}$ is an
orthonormal basis of $P_r$. 
Thus, every $f\in P_r$ can be written uniquely as
$$
f(n)=\sum_{k=1}^{r}\langle f,r^{\scriptscriptstyle {-1/2}}\epsilon_k\rangle
r^{\scriptscriptstyle {-1/2}}\epsilon_k(n),
$$
where 
$$
\langle f, r^{\scriptscriptstyle {-1/2}}\epsilon_k\rangle 
=\sum_{n=1}^{r}f(n)\overline{r^{\scriptscriptstyle {-1/2}}\epsilon_k(n)}
=r^{\scriptscriptstyle {-1/2}}\sum_{n=1}^{r}f(n)\epsilon_k(-n).
$$
This proves (\ref{y:IDFT}) and (\ref{y:DFT}).

\section{DFT and IDFT for even functions $\pmod r$}

Let $E_r$ denote the set of all even functions
$\pmod r$. 
The set $E_r$ forms a complex vector space under the usual addition  
and scalar multiplication. 
In fact, $E_r$ is a subspace of $P_r$. 
Thus (\ref{y:IDFT}) and (\ref{y:DFT}) hold for $f\in E_r$.  
We can also present (\ref{y:IDFT}) and (\ref{y:DFT}) for $f\in E_r$ in 
terms of Ramanujan's sum as is shown below. 

Note that Ramanujan's sum $C(n, r)$ is an integer for all $n$ and can be evaluated by
addition and  subtraction of integers. 
In fact, $C(n, r)$ can be written as $C(n, r)=\sum_{d|(n, r)} d \mu(r/d)$, 
where $\mu$ is the M\"obius function.

An arithmetical function $f\in E_r$ is completely determined by its values 
$f(d)$ with $d|r$. 
Thus $E_r$ is isomorphic to $\C^{\tau(r)}$, where $\tau(r)$ is the number of 
divisors of $r$. 
The inner product (\ref{y:innerp}) in $P_r$ can be written in $E_r$ 
in terms of the Dirichlet convolution. 
In fact, we have 
\begin{equation}\label{y:innerp_phi}
\sum_{k=1\atop (k, r)=d}^{r} 1 = \sum_{j=1\atop (j, r/d)=1}^{r/d} 1 = \phi(r/d), 
\end{equation}
where $\phi$ is Euler's totient function, 
and thus (\ref{y:innerp}) can be written for $f, g\in E_r$ as 
$$
\langle f, g\rangle = \sum_{k=1}^r f(k) \overline{g(k)}
= \sum_{d|r} f(d) \overline{g(d)}\phi(r/d) = (f\overline{g}\ast\phi)(r),  
$$
where $\ast$ is the Dirichlet convolution. 

\begin{theorem}\label{t:orto}
The set
\begin{equation}\label{y:kanta}
\{ (r\phi(d))^{-{1\over 2}}C(\cdot,d)\colon d\mid r\}
\end{equation}
is an orthonormal basis of the inner product space $E_r$.
\end{theorem}

{\sl Proof}\  \ As the dimension of the inner product space $E_r$ is $\tau(r)$
and the number of elements in the set (\ref{y:kanta}) is $\tau(r)$, 
it suffices to show the set (\ref{y:kanta}) is an orthonormal subset of $E_r$.
This follows easily from the relation
$$
\sum_{e\mid r}C(r/e,d_1)C(r/e,d_2)\phi(e)=
 \left\{
   \begin{array}{ll}
    r\phi(d_1)  & \mbox{if $d_1=d_2$,} \\
    0  & \mbox{otherwise,}
   \end{array}
 \right.
$$
where $d_1\mid r$ ja $d_2\mid r$  (see \cite[p. 79]{M}). 
$\Box$

\medskip 

We now present (\ref{y:IDFT}) and (\ref{y:DFT}) for $f\in E_r$. 

\begin{theorem}\label{t:DRT} 
Every $f\in E_r$ can be written uniquely as 
\begin{equation}\label{y:IDRT}
f(n)=r^{-1}\sum_{d\mid r}R_f(d)C(n,d),
\end{equation}
where  
\begin{equation}\label{y:DRT2}
R_f(d)=\phi(d)^{-1}\sum_{n=1}^{r} f(n)C(n, d).
\end{equation}
\end{theorem}

{\sl Proof}\  \ 
On the basis of Theorem \ref{t:orto}, 
\begin{equation}\label{y:f(n)}
f(n)=\sum_{d\mid r}\langle f,(r\phi(d))^{-{1\over 2}}C(\cdot ,d)
\rangle(r\phi(d))^{-{1\over 2}}C(n,d).  
\end{equation}
Applying (\ref{y:innerp}) to (\ref{y:f(n)}) we obtain (\ref{y:IDRT}) 
and (\ref{y:DRT2}). $\Box$

\medskip

The function $R_f$ in (\ref{y:DRT2}) may be referred to as the Discrete
Ramanujan-Fourier Transform of $f$, and (\ref{y:IDRT}) 
may be referred to as the Inverse Discrete Ramanujan-Fourier Transform. 
Cf. \cite{M}. 

Another expression of (\ref{y:DRT2}) can be obtained easily. 
Namely, applying (\ref{y:innerp_phi}) to (\ref{y:f(n)}) and then applying 
$$
\phi(e)C(r/e,d)=\phi(d)C(r/d,e)
$$
(see \cite[p. 93]{M}) we obtain  
\begin{equation}\label{y:DRT22}
R_f(d)=\sum_{e\mid r}f(r/e)C(r/d, e).
\end{equation}
Note that (\ref{y:IDRT}) can also be derived from (\ref{y:IDFT}).  
In fact, if $f\in E_r$, then (\ref{y:DFT}) can be written as 
\begin{eqnarray*}
F_f(k)&=&\sum_{n=1}^{r}f(n)\exp(-2\pi ikn/r)\\
&=&\sum_{e\mid r} \sum_{n=1\atop (n, r)=e}^{r} f(e) \exp(-2\pi ikn/r)\\
&=&\sum_{e\mid r} f(e) \sum_{m=1\atop (m, r/e)=1}^{r/e} \exp(-2\pi ikm/(r/e))\\
&=&\sum_{e\mid r} f(e) C(k, r/e). 
\end{eqnarray*}
A similar argument shows (\ref{y:IDRT}) with $R_f(d)=F_f(r/d)$. 
We omit the details.

\section{The Cauchy product}

It is well known that if $h$ is the Cauchy product of $f\in P_r$ and $g\in P_r$, 
then $F_h=F_f F_g$. 
This follows from the property
$$
\sum_{a+b\equiv n \ppmod{r}} \epsilon_k(a)\epsilon_j(b)=
 \left\{
   \begin{array}{ll}
    r \epsilon_k(n)  & \mbox{if $k\equiv j \ppmod{r}$,} \\
    0  & \mbox{otherwise.}
   \end{array}
 \right.
$$
Analogously, if $h$ is the Cauchy product of $f\in E_r$ and $g\in E_r$, 
then $R_h=R_f R_g$. 
This follows from the property
$$
\sum_{a+b\equiv n \ppmod{r}} C(a, d_1)C(b, d_2)=
 \left\{
   \begin{array}{ll}
    r C(a, d_1)  & \mbox{if $d_1=d_2$,} \\
    0  & \mbox{otherwise,}
   \end{array}
 \right.
$$
where $d_1\mid r$ ja $d_2\mid r$  (see \cite[p. 333]{Si}).


\end{document}